\documentclass[12pt,letterpaper]{article}
\usepackage{amsmath,amsthm}

\font\bbb=msbm10 scaled 1200
\font\bbbs=msbm6 scaled 1000

\DeclareFontShape{OT1}{cmr}{bx}{sc}
      {
      <-> cmbxsc10
      }{}

\newcommand{\bb}[1]{\mbox{\bbb #1}}
\newcommand{\bbs}[1]{\mbox{\bbbs #1}}

\newcommand{\ideal}[1]{\hspace{-1pt} \mbox{ } _{\rm id}  \langle #1 \rangle}

\theoremstyle{plain}
\newtheorem{thm}{Theorem}[section]
\newtheorem{lem}[thm]{Lemma}
\newtheorem{prop}[thm]{Proposition}
\newtheorem{cor}[thm]{Corollary}

\newtheorem{conj}[thm]{Conjecture}

\newtheorem{tetel}{Theorem}

\theoremstyle{remark}

\begin{document}

\title{Applications of Commutator-Type Operators to 
$p$-Groups\thanks{2000 Mathematics Subject 
Classification 20D, 20F, 17B}}

\author{G\'abor Luk\'acs\thanks{I wish to thank the anonymous donors
whose generosity towards the Technion enabled me to pursue my
Master's studies and to do this research.  I am also grateful to York
University and AVI Fellowships for financial support.}
}

\maketitle

\begin{abstract}
For a $p$-group $G$ admitting an automorphism $\varphi$ of
order $p^n$ with exactly $p^m$ fixed points such that
$\varphi^{p^{n-1}}$ has exactly $p^k$ fixed points, we prove that
$G$ has a fully-invariant subgroup of $m$-bounded nilpotency class
with $(p,n,m,k)$-bounded index in $G$. We also establish its analogue for
Lie $p$-rings. The proofs make use of the theory of commutator-type
operators.
\end{abstract}

\section{Introduction}

In this paper all groups and rings are finite, and $p$ is always a
prime number.

\subsection{Brief Historical Background}

In \cite{Ko} \& \cite{Kh3} E. I. Khukhro formulated the following
two conjectures on the structure of $p$-groups admitting an
automorphism of $p$-power order:

\begin{conj}\label{c1}A finite $p$-group admitting an automorphism
of order $p$ with $p^m$ fixed points contains a subgroup of
$m$-bounded class with $(p,m)$-bounded index.
\end{conj}

\begin{conj}\label{c2}A finite $p$-group admitting an automorphism of
order $p^n$ with $p^m$ fixed points contains a subgroup of $m$-bounded
derived length with\linebreak $(p,n,m)$-bounded index.
\end{conj}

(We use the term ``$(a,b,\ldots)$-bounded" for
``bounded above by some function of $a,b,\ldots$".)
R. Shepherd \cite{She}, C. R. Leedham-Green and S. McKay's
\cite{LGM} result gave a positive answer to Conjecture \ref{c1}
when $m=1$. I. Kiming \cite{Ki} gave a positive answer to
Conjecture~\ref{c2} in the case $m=1$. He proved that in this case
there also exists a subgroup of class at most 2 with
$(p,n)$-bounded index.

\bigskip
In \cite{med-1} Yu. Medvedev made two conjectures, on Lie rings:

\begin{conj}\label{c3}A Lie $p$-ring admitting an automorphism of
order $p$ with $p^m$ fixed points contains a nilpotent subring of
m-bounded class whose index is (p,m)-bounded.
\end{conj}

\begin{conj}\label{c4}A Lie $p$-ring admitting an automorphism of
order $p^n$ with $p^m$ fixed points contains a solvable subring of
m-bounded derived length whose index is (p,n,m)-bounded.
\end{conj}

\noindent Yu. Medvedev \cite{med-1} proved that
Conjecture~\ref{c3} implies Conjecture~\ref{c1} and 
Conjecture~\ref{c4} implies Conjecture~\ref{c2}. In \cite{med-2}
Yu. Medvedev deduced Conjecture~\ref{c3}, hence gave a positive
answer to Conjecture~\ref{c1}. In \cite{uj} A. Jaikin-Zapirain proved
Conjecture~\ref{c4}, and therefore gave a positive answer to
Conjecture~\ref{c2}. For a wider historical background I refer to the
Introduction of \cite{med-1} (and \cite{cikk1}).

\subsection{On This Paper}

As we mentioned in \cite{cikk1}, the main results of this paper   
are the following two theorems inspired by Conjecture~\ref{c4} and 
Conjecture~\ref{c2}:

\begin{tetel}\label{app:tetel:main}
Suppose that $L$ is a Lie $p$-ring admitting an automorphism
$\varphi$ of order $p^n$ with exactly $p^m$ fixed points, such
that $\varphi^{p^{n-1}}$ has exactly $p^k$ fixed points. Then $L$
has a nilpotent fully-invariant ideal of $m$-bounded class which
has $(p,n,m,k)$-bounded index in $L$.
\end{tetel}

\begin{tetel}\label{group:tetel}
If a  finite $p$-group $P$ admits an automorphism of order $p^n$ with
exactly $p^m$ fixed points, such that $\varphi^{p^{n-1}}$ has
exactly $p^k$ fixed points, then $P$ has a fully-invariant
subgroup of $(p,n,m,k)$-bounded index which is nilpotent of
$m$-bounded class.
\end{tetel}

These theorems give in some sense a stronger result ($m$-bounded
nilpotency class) than Conjectures~\ref{c2} and~\ref{c4}, however we have
to ``pay" for that by involving an additional parameter $k$ into the bound
for the index. We always have $m \leq k$. Note that if $m=k$,
Theorems~\ref{app:tetel:main} and~\ref{group:tetel} are obvious
consequences  of Medvedev's theorems mentioned above. But in the case
where $m < k$, we obtain a much better bound for the class ($m$-bounded).

\bigskip
In Section~\ref{SECT:APP} we prove Theorem~\ref{app:tetel:main}. The
proof makes use of the framework of commutator-type operators,   
which we developed in \cite{cikk1}, especially Theorem~III (to
which we will refer in this paper as ``Theorem~III"). It follows 
the design of Khukhro's version \cite[14.2]{konyv} of a proof of a
theorem of Medvedev (Conjecture~\ref{c3}) mentioned above. 

In Section~\ref{SECT:GRP} we prove that Theorem~\ref{group:tetel}
follows from Theorem~\ref{app:tetel:main}. The proof is done in a
similar way as suggested by the referee in \cite{med-1}.

\bigskip The experienced reader may wish to omit Section 2 (except for
notation and definitions) and jump directly to Section 3.

\section{Preliminary Facts}
For the convenience of the reader and further reference we mention
some results and notations which will be used occasionally in the
paper.

\subsection{Automorphisms of Abelian Groups}

For any positive integer $k$ we denote by $\phi_k(x)$ the
cyclotomic polynomial of order $p^k$.

Let $G$ be a group and $\varphi$ an automorphism of the group.
The set of fixed points of $\varphi$ on $G$ is $C_G(\varphi)$.

Let $A$ be a finite additive $p$-group, that admits an
automorphism $\varphi$ of order $p^n$ with exactly $p^m$ fixed
points.

\begin{prop}{\rm\cite[3.2]{thesis}}
\label{prep:prop:eq:cyclo}
We have
\begin{equation}
\label{prep:eq:cyclo}
a+a \varphi + a \varphi^2 + \cdots + a \varphi^{p^n-1}
= a \prod\limits_{k=1}^n \phi_k(x) = 0
\end{equation}
for all $a \in p^m A$.
\end{prop}

\begin{lem}{\rm\cite[3.2]{cikk1}}\label{prep:cyc:fix} Suppose that
$\varphi$ satisfies (\ref{prep:eq:cyclo}) on $A$. Then:

\parskip 0pt

{\rm (a)}
for any $\varphi$-invariant section $U$ of $A$ we
have $p^nC_U(\varphi)=0$;

\hangindent 37pt
{\rm (b)}
for any homocyclic $\varphi$-invariant section $V$ of $A$
of exponent $p^s$
we have   
$|C_V(\varphi)| = |C_{p^i V / p^{i+n}V} (\varphi) |$ whenever
$0 \leq i \leq s - n$.
\end{lem}

\begin{lem}{\rm\cite[3.4]{cikk1}}\label{prep:lemma:n'} Let $A$ be
an abelian group. If $\phi_{n^\prime}(\varphi)=0$ and $\varphi$
has order $p^n$ on $A$ then either $n^\prime=n$ or $n^\prime > n$.
If $n^\prime > n$, then $pA =0$.
\end{lem}

\subsection{Lie Rings}

\begin{prop}\label{basic:prop:p+1}{\rm\cite[7.20]{konyv}}
Let $L$ be a Lie ring of derived length 2, admitting an
automorphism $\varphi$ of order $p$. Then
$\gamma_{p+1}(pL) \subseteq \! \! \ideal{C_L (\varphi)}$.
\end{prop}

A Lie ring is called Lie $p$-ring if its additive structure is a
$p$-group.

\section{\label{SECT:APP}Proof of Theorem~\ref{app:tetel:main}}

%In this section we will use the tool of commutator-type operators
%and the framework we developed in \cite{cikk1} to prove the
%following theorem:

The proof of Theorem~\ref{app:tetel:main} follows the design of
\cite[14.2]{konyv}: we split it up into two theorems
(Theorem~\ref{app:theorem:deriv} and  Theorem~\ref{app:theorem:class}
below). Then, in subsection~\ref{app:subsect:proof} we show
that they imply Theorem~\ref{app:tetel:main}. It is interesting to note 
that the tool of commutator-type operators we developed in \cite{cikk1}
gives a common framework for proving Theorem~\ref{app:theorem:deriv} and
Theorem~\ref{app:theorem:class}.

\subsection{\label{app:subsect:deriv}$m$-Bounded Derived Length}

\begin{thm}\label{app:theorem:deriv}
%{\rm\cite[14.3]{konyv}**}
Suppose that $L$ is a Lie $p$-ring. If $L$ admits an automorphism
$\varphi$ of order $p^n$ with exactly $p^m$ fixed points, such
that $\varphi^{p^{n-1}}$ has exactly $p^k$ fixed points, then $L$
has a soluble fully-invariant ideal of $m$-bounded derived length
which has $(p,n,m,k)$-bounded index in $L$.
\end{thm}

\begin{proof}$\varphi^{p^{n-1}}$ has $p^k$ fixed points, thus by
Proposition~\ref{prep:prop:eq:cyclo} $\varphi^{p^{n-1}}$ satisfies
the polynomial $\phi_1(x)=0$ on $p^k L$ (for its order is $p$).
This means that $\phi_1(\varphi^{p^{n-1}}) = \phi_n(\varphi)=0$ on
$p^k L$. $p^k L$ has $(p,n,m,k)$-bounded index in $L$ (because $L$ has
$(p,n,m)$-bounded rank, see \cite[2.7]{konyv}), hence replacing $L$
by $p^k L$ we may assume that $\phi_n(\varphi)=0$ on $L$. We note that
this is the only step where $k$ plays any role. Later we will have
$(p,n,m)$-bounded index.

If $R$ is a Lie-ring, clearly $T(A_1,A_2)=[A_1,A_2]$ is a
commutator-type operator in $2$ variables on $R$.

Let ${\cal L}$ be the set of all the triplets consisting of a
finite Lie $p$-ring and the Lie bracket $[\cdot,\cdot]$ as a
commutator-type operation in $2$ variables and an automorphism
$\psi$ of $L$ of $p$-power order. 

We show that ${\cal L}$ satisfies condition (c) of Theorem~III
(clearly ${\cal L}$ is closed under "taking" $(T,\psi)$-invariant
sections). Let $(R,T,\varphi) \in {\cal L}$ such that $R$ is
homocyclic 
and admits an automorphism $\varphi$ of order $p^n$ with $p^m$ fixed
points, such that $\phi_n(\varphi)=0$.
Define $[[ \cdot, \cdot]]$ as in \cite[2.14]{thesis} (or \cite[13.24]{konyv})
(with $t$ in place of $s$). Let $\widetilde R$ be the ``lifted" ring, and
let $\bb{T}$ be the ``lifted operator" that 
$T$ defines on $\widetilde R$ (see \cite[page 15]{cikk1}). 
For $A,B \leq R$, clearly
$\bb{T} (\widetilde A,\widetilde B) = [[\widetilde A,\widetilde B]]$,
because by the definion of $[[ \cdot, \cdot]]$ we have
$p^t [[A,B]] = [A,B]=T(A,B)$.
Hence $\gamma_u(\widetilde R) = \gamma_u^{\bbs{T}}(\widetilde R)$.
%%(Note that in this case $l=2$.)
By Proposition~\cite[13.26]{konyv}, $\widetilde R$ forms a
Lie-ring with this new multiplication. (All these are true, because the  
``classical" top of $R$ used in~\cite{konyv} 
coincides with $t(R)$ that $T$ defines.)

$\phi_n(\varphi)=0$ on $R$, thus $\phi_n(\varphi)=0$ on
$\widetilde R$ (since it is a $\varphi$-invariant section of~$R$).
$pC_{\widetilde R}(\varphi^{p^{n-1}})=0$, by
Lemma~\ref{prep:cyc:fix} applied to $\varphi^{p^{n-1}}$ (because
$\varphi^{p^{n-1}}$ satisfies $\phi_1(x)=0$ on $\widetilde R$).
Let $h=h(p)$ (the Higman number). By Higman's theorem, we have   
\vspace{-5pt}
\[
p^{h+2} \gamma_{h+1}(\widetilde R) = p
\gamma_{h+1}(p \widetilde R) \subseteq p \ideal{C_{\widetilde R}
(\varphi^{p^{n-1}})} =\ideal{p C_{\widetilde R}
(\varphi^{p^{n-1}})}=0 \mbox{.}
\]
Thus, for $u(p,n,m)=h(p)$,
$v(p,n,m)=h(p)+2$, ${\cal L}$ satisfies the assumptions of
Theorem~III. Note, that in this step we made use of the fact that
$\widetilde R$ is also a Lie-ring (otherwise we could not apply
Higman's theorem).

Therefore, by Theorem~III, there exist  an $m$-bounded number $g=g(m)$
and a $(p,n,m)$-bounded number $r=r(p,n,m)$, such that for any Lie
$p$-ring $L$ admitting an automorphism $\varphi$ of order $p^n$ with $p^m$
fixed points such that $\phi_n(\varphi)=0$, we have
$(p^r L)^{(g)} = T^{(g)}(p^r L)=0$. $p^rL$ has $(p,n,m)$-bounded
index in $L$ (because the rank of $L$ and $r$ are $(p,n,m)$-bounded).

\bigskip\noindent This completes the proof of
Theorem~\ref{app:theorem:deriv}.
\end{proof}

\subsection{\label{app:subsect:class}$m$-Bounded Nilpotency Class}

To shorten the notation, we will write
$[C,_i D] = [C, \underbrace{D,\ldots,D}_{i\mbox{ times}}]$.

\begin{thm}\label{app:theorem:class}
Let $L = A \oplus B$ be a Lie $p$-ring, with abelian ideal $A$ and
abelian subgroup $B$. Suppose that $L$ admits an automorphism
$\varphi$ of order $p^n$ with exactly $p^m$ fixed points, such
that $A$ and $B$ are $\varphi$-invariant and
$\varphi^{p^{n-1}}$\hspace{-4pt} has exactly $p^k$\hspace{-1pt}
fixed points on $L$.\hspace{-2.5pt} Then there exits a
$(p,n,m,k)$-bounded number $f=f(p,n,m,k)$ and an $m$-bounded number
$g=g(m)$ such that\linebreak $p^{f}[A,_g B]=0$.
\end{thm}

\begin{proof}Since the order of $\varphi^{p^{n-1}}$ is $p$, by
Proposition~\ref{prep:prop:eq:cyclo} applied on
$\varphi^{p^{n-1}}$ we obtain that $\phi_1(\varphi^{p^{n-1}})=
\phi_n(\varphi) = 0$ on $p^k L$.
Since
\begin{equation} \label{app:eq:pc}
[p^k A,_{g(m)} (p^k B)]=
p^{(g(m)+1)k} [A,_{g(m)}B]\mbox{,}
\end{equation}
replacing $L$ by $p^k L$ we may assume that $\phi_n(\varphi)=0$ on
$L$ from the outset. We note that this is the only step where $k$
plays any role. Later we will have $(p,n,m)$-bounded index.

Let ${\cal A}$ be the set of the triplets $(A,T,\psi)$ such that
$T$ is the operator $T(C)=[C,B]$ on $A$, where $[\cdot,\cdot]$ is
the Lie-bracket in a Lie $p$-ring $L=A\oplus B$ (so $A$ and $B$ 
are $p$-groups), $A$ is an abelian $\psi$-invariant ideal in $L$,
$B$ is an abelian $\psi$-invariant subring of $L$ and $\psi$ is an
automorphism of $L$ such that $\phi_{n^\prime}(\psi)=0$ on $L$ for
some $n^\prime$. From $\phi_{n^\prime}(\psi)=0$ it follows that
$\psi^{p^{n^\prime}}=1$, so $\psi$ has $p$-power order.

We claim that ${\cal A}$ satisfies the conditions of Theorem~III.
First, suppose that $(A,T,\psi) \in {\cal A}$ and $S=C/D$ is a
$(T,\psi)$-invariant section of $A$. Then $S$ is "$B$-invariant"
(i.e. $[S,B] \leq S$), thus $C$ and $D$ are $B$-invariant, hence
$D$ is an ideal of $C \oplus B$. So $S \oplus B$ forms a new
Lie-ring, being a Lie-section of $L$ (a subring modulo its
ideal). Clearly $(S,T,\psi) \in {\cal A}$.

Let $(A,T,\varphi) \in {\cal A}$ such that $A$ is homocyclic,
$|C_A(\varphi)| = p^m$ and the order of $\varphi$ on $L$ is $p^n$.
Let $\widetilde L = \widetilde A \oplus B$ with the operation
$[[\cdot,\cdot]]$: for any $x \in A$ and $y \in B$ we put
$[[x,y]]$ to be a $p^{t}$th root of $[x,y]$, where $t=t(A)$. If either
$z_1,z_2 \in A$ or $z_1, z_2 \in B$ we define $[[z_1,z_2]]=0$. 
The additive factor group $\widetilde{L}$ endowed with this operation
($[[\widetilde{x},\widetilde{y}]]=\widetilde{[[x,y]]}$, where
tilde denotes image in $\widetilde{L}$), is a Lie ring, and the
automorphism of the additive group of $L$ induced by $\varphi$ is
an automorphism of the Lie ring $\widetilde{L}$ 
(see Proposition~\cite[6.7, page 41]{thesis}).

We show now that condition (c) of Theorem~III holds here. Since  
$\widetilde A$ and $B$ are abelian and $\widetilde A$ is an ideal 
in $\widetilde L$, we have
$\gamma_{u+1}(\widetilde L) = [[A,_u B ]]$.
Clearly $\bb{T}(C) = [[C, B]]$  for any $C \leq \widetilde A$, thus
\begin{equation} \label{app:eq:1s}
\gamma_{u+1}(\widetilde L) = \gamma_{u+1}^{\bbs{T}} (\widetilde A)
\mbox{.}
\end{equation}

By the definition of ${\cal A}$, $\phi_{n^\prime} (\varphi) = 0$
on $L$. By Lemma~\ref{prep:lemma:n'}, if $n \neq n^\prime$, then
$pL=0$, so we are done. Thus we may assume that $n=n^\prime$. As
an abelian group, $\widetilde L$ is a $\varphi$-invariant section
of $L$, hence $\phi_n (\varphi) = 0$ on $\widetilde L$. Let
$\sigma =\varphi^{p^{n -1}}$ (as an automorphism of $\widetilde
L$). Then $\phi_1 (\sigma)=0$, thus (applying
Lemma~\ref{prep:cyc:fix}(a) on $\sigma$)
\begin{equation} \label{app:eq:2s}
p C_{\widetilde L}(\sigma) = 0 \mbox{.}
\end{equation}

The ``lifted" Lie-ring $\widetilde L$ has derived length $2$,
hence by Proposition~\ref{basic:prop:p+1} we have
$\gamma_{p+1}(p\widetilde L) \subseteq \! \!
\ideal{C_{\widetilde L} (\sigma)}$.
Hence,
\[
p^{p+2} \gamma_{p+1}^{\bbs{T}} (\widetilde A)
\stackrel{(\ref{app:eq:1s})}{=} p^{p+2} \gamma_{p+1}
(\widetilde L)
=  p \gamma_{p+1}(p\widetilde L)
\stackrel{(\ref{basic:prop:p+1})}{\subseteq} 
p\ideal{C_{\widetilde L} (\sigma)}
= \!\ideal{p C_{\widetilde L}(\sigma)}
\stackrel{(\ref{app:eq:2s})}{=} 0 \mbox{.}
\]
Therefore, ${\cal A}$ satisfies the conditions of
Theorem~III with $v(p,n,m) = p+2$, $u(p,n,m)=p$.

So there exist a $(p,n,m)$-bounded number $f(p,n,m)$ and an
$m$-bounded number $g(m)$ such that for any $L=A \oplus B$ such
that $\phi_n (\varphi) = 0$ on $L$ and $\varphi$ has at most $p^m$
fixed points {\bf on $\boldsymbol{A}$}, then $T^{(g)}(p^f A) = p^f
T^{(g)} (A)$ (in particular, it is true if $|C_L(\varphi)| = p^m$).

\bigskip\noindent
This completes the proof of Theorem~\ref{app:theorem:class}.
\end{proof}

%\subsection{\label{app:subsect:proof}Proof of Theorem~\ref{app:tetel:main}}
\subsection{\label{app:subsect:proof}
Theorem \ref{app:theorem:deriv} and \ref{app:theorem:class}
imply Theorem~\ref{app:tetel:main}}

In this subsection we will show, how Theorem~\ref{app:tetel:main}
follows from Theorem~\ref{app:theorem:deriv} and
Theorem~\ref{app:theorem:class}. 
First note the following corollary of Theorem~\ref{app:theorem:class}.

\begin{cor}\label{app:cor:class-deriv}
%{\rm\cite[14.39]{konyv}}
Let $L$ be a soluble Lie $p$-ring of derived length $2$.  Suppose that
$L$ admits an automorphism $\varphi$ of order $p^t$ with exactly
$p^s$ fixed points, such that $\varphi^{p^{t-1}}$ has $p^r$ fixed
points. Then $\gamma_{g(s)} ( p^{f(p,t,s,r)} L) =0$ for a
$(p,t,s,r)$-bounded number $f(p,t,s,r)$ and an $s$-bounded number
$g(s)$.
\end{cor}

This corollary is an adaption of \cite[14.39]{konyv} to our conditions.

\begin{proof}Put $A = [L,L]$ and $B = L / [L,L]$. We define the (Lie)
operation of $A \oplus B$ as follows:  if both $x$, $y$ in $A$ or
in $B$, then $[x,y]=0$, and if $y \in B$ then $[a,b] = [a,y]$ if
$b = y + [L,L] \in B$. $A$ and $B$ are $\varphi$-invariant, and we
have $|C_{A \oplus B}(\varphi)|=p^m \leq p^{2s}$, and $|C_{A
\oplus B}(\varphi^{p^{t-1}})|=p^k \leq p^{2r}$ (see \cite[2.12]{konyv}). 
It follows form the definition of
$L$, that $A$ is an abelian ideal in $A \oplus B$ and $B$ is an   
abelian subring, and they are $\varphi$-invariant.

If $\varphi$ has order less than $p^t$ on $A \oplus B$, then
$\varphi^{p^{t-1}}=1$ on $A \oplus B$, so we have $p^{2r} [L,L]=0$
(for $|C_{A \oplus B}(\varphi^{p^{t-1}})|\leq p^{2r}$), and we are
done. So we may assume that $\varphi$ has order $p^t$ on $A \oplus
B$. Thus $A \oplus B$ satisfies the conditions of
Theorem~\ref{app:theorem:class}, with parameters $t$, $m \leq 2s$
and $k \leq 2r$.
Hence for some $s$-bounded number $v$ and $(p,t,s,r)$-bounded number $u$
we have $ p^{u(p,t,s,r)} [A,_{v(s)} B]=0$.
($m$ and $k$ appear in the bounds given by
Theorem~\ref{app:theorem:class}, but they are $s$ and $r$ bounded 
respectively and the bounds can be assumed to be  monotonic functions). By
the definition of the operation in $A \oplus B$, it implies that
\[
\gamma_{v + 1}(p^{[u/(v+1)] +1} L ) \leq
p^u [ \hspace{2pt} [L,L],_{v(s)}L ]  =
p^{u(p,t,s,t)} [A ,_{v(s)} B]=0\mbox{.}
\]
This completes the proof of Corollary~\ref{app:cor:class-deriv}. 
\end{proof}

\begin{prop}\label{app:prop:class-deriv}
%{\rm\cite[14.40]{konyv}}
Suppose that $L$ is a soluble Lie $p$-ring of derived length $d$.
If $L$ admits an automorphism $\varphi$ of order $p^n$ with
exactly $p^m$ fixed points, such that $\varphi^{p^{n-1}}$ has
$p^k$ fixed points on $L$, then $L$ has a fully-invariant
nilpotent ideal of $(m,d)$-bounded class which has
$(p,n,m,k,d)$-bounded index in $L$.
\end{prop}

This proposition is an adaption of \cite[14.40]{konyv} to our conditions.

\begin{proof}We proceed by induction on $d$, the derived length of
$L$. If $d=2$, the result follows from
Corollary~\ref{app:cor:class-deriv}. Suppose that $d > 2$. By the
inductive hypothesis, $\gamma_v ( p^u [L,L]) = 0$ for
$u=u(p,n,m,k,d-1)$ and $v = v(m,d-1)$.  We have
$|C_{[L,L]}(\varphi^{p^{n-1}})|, |C_{L/L^{(2)}}(\varphi^{p^{n-1}})|
\leq |C_L(\varphi^{p^{n-1}})|$ (see \cite[2.12]{konyv}), so
by  Corollary~\ref{app:cor:class-deriv} (applied to $L/ L^{(2)}$),
$\gamma_g(p^f L) \leq L^{(2)}$ for some  $g = g(m)$ and
$f = f(p,n,m,k)$. We put $w = \max \{[u/2]+1 ,f \}$, which   
is $(p,n,m,k,d)$-bounded, and put $M = p^w L$. Then
\begin{eqnarray*}
\gamma_{g+4} (M) &  = & \gamma_{g + 4} (p^w L) \leq
[L^{(2)},p^w L, p^w L, p^w L, p^w L] \\
& \leq & p^{4w} L^{(2)} = (p^w L)^{(2)} = M^{(2)}\mbox{,}
\end{eqnarray*}
therefore $\gamma_{g+4} (M)  \leq M^{(2)}$. 

But since $u \leq 2w$,
$\gamma_v([M,M])=\gamma_v (p^{2w}[L,L]) \leq\gamma_v (p^u [L,L]) =0$.
Applying \cite[5.27]{konyv} with $L = M$ and $N=[M,M]$ we obtain that the
nilpotency class of $M$ is bounded in the terms of $v$ and $g$, that is
$(m,d)$-bounded. Thus, $M$ is the required fully-invariant ideal of
$(m,d)$-bounded class with $(p,n,m,k,d)$-bounded index in $L$.
\end{proof}

\begin{proof}[\bf\scshape Proof of Theorem~\ref{app:tetel:main}.]According
to Theorem~\ref{app:theorem:deriv}, $L$ has a soluble fully-invariant
ideal of $m$-bounded derived length with $(p,n,m,k)$-bounded index
in $L$. Thus, without loss of generality we may assume that $L$
has $m$-bounded derived length.  By
Proposition~\ref{app:prop:class-deriv}, $L$ has a nilpotent 
fully-invariant ideal of $(m,d)$-bounded class with
$(p,n,m,k,d)$-bounded index in $L$. Since $d$ is $m$-bounded, the
result follows.

\bigskip\noindent
This completes the proof of Theorem~\ref{app:theorem:deriv}.
\end{proof}

\section{\label{SECT:GRP}Proof of Theorem~\ref{group:tetel}}

The structure of the proof of Theorem~\ref{group:tetel}
follows the scheme suggested by the referee in~\cite{med-1}.

\begin{lem}\label{group:lemma}If
$G$ satisfies the conditions of the Theorem, then there exist a
$(p,n,m,k)$-bounded number $f_1=f_1(p,n,m,k)$ and an $m$-bounded
number $g=g(m)$ such that
$\gamma_g(G^{p^{f_1}}) \subseteq \gamma_p (G)$.
\end{lem}

\begin{proof}First suppose that, $G$ is of class $ < p$. Let $L$ be
the corresponding Lie-ring, by the Lazard correspondence
\cite[chapter 10, 10.24]{konyv} (inverse Baker-Hausdorff formula).
$L$ fulfils the conditions of Theorem~\ref{app:tetel:main}, so we
have $\gamma_{g(m)}(p^{f_1(p,n,m,k)}L)=0$ for the
$(p,n,m,k)$-bounded number $f_1$ and the $m$-bounded number $g$.
The correspondent of $p^{f_1(p,n,m,k)}L$ is
$G^{p^{f_1(p,n,m,k)}}$, so $\gamma_g(G^{p^{f_1}})=1$.

Now suppose that $G$ is an arbitrary group satisfying the
conditions of the theorem. Then, applying the consideration above
to $\overline G = G/ \gamma_p(G)$ we obtain that
$\gamma_g(G^{p^{f_1}})\gamma_p(G)/\gamma_p(G) =
\gamma_g({\overline G}^{p^{f_1}}) = 1$.
Therefore $\gamma_g(G^{p^{f_1}}) \subseteq \gamma_p (G)$, as desired.
\end{proof}

\begin{proof}[\bf\scshape Proof of Theorem~\ref{group:tetel}.]Let 
$\psi = \varphi^{p^{n-1}}$, then $|C_P (\psi)| = p^k$.  $P$ has a fully
invariant subgroup of $(p,k)$-bounded index which is nilpotent of class at
most $h(p)$ (see \cite[8.1]{konyv}; it follows from the proof that the
characteristic subgroup found in the theorem is, in fact, 
fully-invariant). Thus, we may assume that
\begin{equation} \label{group:assume:h(p)}
\mbox{$P$ is nilpotent of class at most $h(p)$.}
\end{equation}

Let $g=g(m)$ and $f_1=f_1(p,n,m,k)$ from Lemma~\ref{group:lemma} applied
to $P$. We show that there exists a $(p,n,m,k)$-bounded number
$f=f(p,n,m,k)$, such that $\gamma_g (P^{p^f}) = 1$. If   
$p \leq g(m)$, then the class of $P$ is $m$-bounded  
(by~(\ref{group:assume:h(p)})). If $p > g(m)$, we proceed by 
induction on $c$, the class of $P$. By Lemma~\ref{group:lemma},  
$\gamma_g(P^{p^{f_1}}) \subseteq \gamma_p(P)$. So
$\gamma_{c-p+g}(P^{p^{f_1}}) \subseteq \gamma_c(P) =1$. Since
$g < p$, $c-p+g < c$, thus by the inductive hypothesis there
exists a $(p,n,m,k)$-bounded number $f_2(p,n,m,k)$ such that
$\gamma_g( (P^{p^{f_1}})^{p^{f_2}}) = 1$.

Let $f=f_1 + f_2$. Since  
$P^{p^{f_1 + f_2}} \subseteq  (P^{p^{f_1}})^{p^{f_2}}$, we have
$\gamma_g(P^{p^{f}})=1$. 
Clearly $P^{p^f}$ is fully-invariant. Let $H = P/ P^{p^f}$.
By Burnside Theorem, the number of generators of $H$ is
the same as the number of generators of $K =H/ \Phi(H)$ (here
$\Phi(H)=H^p[H,H]$
is the Frattini subgroup of $H$). $K$ is an
abelian $\varphi$-invariant section of $H$
with $|C_K(\varphi)| \leq p^m$ (see \cite[2.12]{konyv}).
Hence (by Corollary~\cite[2.7]{konyv}), $r(K) \leq mp^n$, so $H$
has at most $m p^n$ generators. The exponent of $H$ is at most
$p^f$. Thus the exponent, the number of generators and the nilpotency
class of $H$ are $(p,n,m,k)$-bounded. Therefore, the order of  $H$ is
$(p,n,m,k)$-bounded (see \cite[6.12(c)]{konyv}), hence $P^{p^f}$ has
$(p,n,m,k)$-bounded index in $P$. By~(\ref{group:assume:h(p)})
$c \leq h(p)$, thus we have $p$-bounded number of steps in the induction.

\bigskip\noindent
This completes the proof of Theorem~\ref{group:tetel}.
\end{proof}

%\newpage
\section*{Acknowledgements}

\mbox{ }

I am deeply indebted to my Master's thesis supervisor,
Prof.~Arye~Juh\'asz, for his dedicated mentorship that made this
research possible.

\bigskip
I am grateful to Prof.~Pavel~Shumyatsky, who gave me preprints of
the papers \cite{med-1} and \cite{med-2} of Medvedev in February
1999, before their publication.

\bigskip
I would like to thank to Prof.~Walter~Tholen, who was of great
assistance in editing the paper -- both in its formation,
articulation, and preparation for submission.

\bigskip
I would like to express my heartfelt gratitude to my father,
Dr.~J\'anos Luk\'acs, for providing me with a secure home
environment making possible the flowering of this research.

\bigskip
Last, but not least, special thanks to my students for their
enormous encouragement.

%%%%%%%%%%%%%%%%%%%%%%%%%%%%%%%%
%%%%%%%%%%%%%%%%%%%%%%%%%%%%%%%%

%% Bibliography

{\bigskip\bigskip\noindent Department of Mathematics \& Statistics\\
York University, 4700 Keele Street\\
Toronto, Ontario, M3J 1P3\\
Canada

\bigskip\noindent{\em e-mail: lukacs@mathstat.yorku.ca} }

\end{document}